\documentclass{article}
\usepackage{amsfonts}
\usepackage{amsmath}

\newtheorem{theorem}{Theorem}

\newtheorem{proposition}[theorem]{Proposition}

\input{tcilatex}
\begin{document}

\title{On the Cheeger-Gromoll metric}
\author{Melek ARAS\thanks{%
Department of Mathematics, Faculty of Arts and Sciences, Giresun Universty,
28049 Giresun,Turkey\textbf{e-mail:}%
melekaras25@hotmail.com;melek.aras@giresun.edu.tr} \\
Giresun,Turkey }
\maketitle

\begin{abstract}
The purpose of this paper is to investigate applications the covariant
derivatives, killing vector fields and to calculate the components of the
curvature tensor $^{CG}R$ of the Cheeger-Gromoll metric with respect to
adapted frames in a the Riemannian manifold to its tangent bundle $T(M_{n})$
.


\textbf{AMS 2010:} 53C22, 53C25
\end{abstract}

\bigskip \textbf{1. Introduction}

Let $\left( M_{n},g\right) $ be a Riemannian manifold and $T\left(
M_{n}\right) $ its tangent bundle with the projection $\pi :T\left(
M_{n}\right) \rightarrow M_{n}\cite{8}$. In the present paper $\Im
_{q}^{p}\left( M_{n}\right) $ is the set of all tensor fields of type $%
\left( p,q\right) $ on $\left( M_{n}\right) $.

Cheeger and Gromoll studied complete manifolds of non-negative curvature in 
\cite{2}. In spired by the paper \cite{2} of Cheeger and Gromoll, Musso and
Tricerri \ defined the metric $^{CG}g$ which they called the $%
Cheeger-Gromoll $ metric on tangent bundle of Riemannian manifold in \cite{4}%
. The Levi-Civita connection of $^{CG}g$ are calculeted by A.A. Salimov \
and S. Kazimova in \cite{5} with respect to the adapted frame. In \cite{1}
M. Abbas and M. Sarih studied killing vector field on tangent bundles with $%
Cheeger-Gromoll$ metric. In \cite{6} Sekizawa calculated Levi-Civita
connection and curvature tensor of the metric $^{CG}g$ (for more detailes
see \cite{3}).

Let there be given a vector field $X=X^{i}\partial _{i}$ and covector field $%
g_{X}=g_{ji}X^{j}dx^{i}$ in $U\subset M_{n}$. Then $\gamma g_{X}\in $ $\Im
_{0}^{0}\left( M_{n}\right) $ is a function on $\pi ^{-1}:\left(
M_{n}\right) \rightarrow T\left( M_{n}\right) $ defined by $\gamma
g_{X}=y^{j}g_{ji}X^{i}$ with respect to the induced coordinates $\left(
x^{i},y^{i}\right) $ (where $\pi $ is the naturel projection $\pi :T\left(
M_{n}\right) \rightarrow M_{n}$)\cite{8}. Now, denote by $r$ the norm a
vector $y=x^{\overline{i}}=y^{i},$i.e., $r^{2}=g_{ji}y^{j}y^{i}$. The metric 
$^{CG}g$ on $T\left( M_{n}\right) $ is given by

\begin{equation*}
\left\{ 
\begin{array}{l}
^{CG}g\left( ^{H}X,^{H}Y\right) =^{V}\left( g\left( X,Y\right) \right) \\ 
^{CG}g\left( ^{H}X,^{V}Y\right) =0, \\ 
^{CG}g\left( ^{V}X,^{V}Y\right) =\frac{1}{\alpha }\left[ ^{V}\left( g\left(
X,Y\right) \right) +\left( \gamma g_{X}\right) \left( \gamma g_{Y}\right) %
\right] ,%
\end{array}%
\right.
\end{equation*}%
for all vector field $X,Y\in \Im _{0}^{1}\left( M_{n}\right) $, where $%
^{V}\left( g\left( X,Y\right) \right) =\left( g\left( X,Y\right) \right)
\circ \pi $ and $\alpha =1+r^{2}.$

Then the special frame is called the adapted frame. The $Cheeger-Gromoll$
metric $^{CG}g$ has components

\begin{equation}
^{CG}g_{\alpha \beta }=\left( 
\begin{array}{cc}
g_{ji} & 0 \\ 
0 & \frac{1}{\alpha }\left( g_{ji}+g_{js}g_{it}y^{s}y^{t}\right)%
\end{array}%
\right)  \label{1}
\end{equation}%
with respect to the adapted frame $\left\{ X_{\left( i\right) },X_{\left( 
\overline{i}\right) }\right\} $. The Levi-Civita connections of the $^{CG}g$ 
\cite{5} are

\begin{equation}
\left\{ 
\begin{array}{c}
^{CG}\Gamma _{ji}^{h}=\Gamma _{ji}^{h},\Gamma _{ji}^{\overline{h}}=-\frac{1}{%
2}R_{jik}^{h}y^{k},^{CG}\Gamma _{j\overline{i}}^{h}=-\frac{1}{2\alpha }%
R_{.jki}^{h.}y^{k} \\ 
^{CG}\Gamma _{j\overline{i}}^{\overline{h}}=\Gamma _{ji}^{h},^{CG}\Gamma _{%
\overline{j}i}^{h}=-\frac{1}{2\alpha }R_{.ikj}^{h.}y^{k},\Gamma _{\overline{j%
}i}^{\overline{h}}=0,\Gamma _{\overline{j}\overline{i}}^{h}=0, \\ 
^{CG}\Gamma _{\overline{j}\overline{i}}^{\overline{h}}=-\frac{1}{\alpha }%
\left( y_{j}\delta _{i}^{h}+y_{i}\delta _{j}^{h}\right) +\frac{1+\alpha }{%
\alpha }g_{ji}y^{h}-\frac{1}{\alpha }y_{j}y_{i}y^{h}%
\end{array}%
\right.
\end{equation}%
with respect to the adapted frame , where $y_{j}=g_{ji}y^{i}$, $%
R_{.ikj}^{h.}=g^{ht}g_{js}R_{tik}^{h}.$

\textbf{2. The metric }$Cheeger-Gromoll$\textbf{\ in adapted frames}

Let $X$ be a vector field in $T\left( M_{n}\right) $ and $\left( X^{\alpha
}\right) =\left( 
\begin{array}{c}
X^{h} \\ 
X^{\overline{h}}%
\end{array}%
\right) $ its components with respect to the adapted frame.Then the
covariant derivative $\nabla X$ has components

\begin{equation}
^{CG}\nabla _{\beta }X^{\alpha }=D_{\beta }X^{\alpha }+^{CG}\Gamma _{\beta
\delta }^{\alpha }X^{\delta },
\end{equation}%
$^{CG}\Gamma _{\beta \delta }^{\alpha }$ being Levi-Civita connections of
the metric $^{CG}g$ with respect to the adapted frame .

Now let us consider the covariant derivatives of vertical lift, complete
lift and horizontal lift. Then we have respectively components

\begin{equation*}
^{CG}\nabla _{\beta }^{V}X^{\alpha }=\left( 
\begin{array}{cc}
^{CG}\nabla _{i}^{V}X^{h} & ^{CG}\nabla _{\overline{i}}^{V}X^{h} \\ 
^{CG}\nabla _{i}^{V}X^{\overline{h}} & ^{CG}\nabla _{\overline{i}}^{V}X^{%
\overline{h}}%
\end{array}%
\right)
\end{equation*}

\begin{equation}
\left\{ 
\begin{array}{l}
^{CG}\nabla _{i}^{V}X^{h}=-\frac{1}{2\alpha }R_{.ikm}^{h.}y^{k}X^{m} \\ 
^{CG}\nabla _{\overline{i}}^{V}X^{h}=0 \\ 
^{CG}\nabla _{i}^{V}X^{\overline{h}}=\nabla _{i}X^{h} \\ 
^{CG}\nabla _{\overline{i}}^{V}X^{\overline{h}}=\left[ -\frac{1}{\alpha }%
\left( y_{i}\delta _{m}^{h}+y_{m}\delta _{i}^{h}\right) +\frac{1+\alpha }{%
\alpha }g_{im}y^{h}-\frac{1}{\alpha }y_{i}y_{m}y^{h}\right] X^{h}%
\end{array}%
\right. ,
\end{equation}

\begin{equation*}
^{CG}\nabla _{\beta }^{C}X^{\alpha }=\left( 
\begin{array}{cc}
^{CG}\nabla _{i}^{C}X^{h} & ^{CG}\nabla _{\overline{i}}^{C}X^{h} \\ 
^{CG}\nabla _{i}^{C}X^{\overline{h}} & ^{CG}\nabla _{\overline{i}}^{C}X^{%
\overline{h}}%
\end{array}%
\right) ,
\end{equation*}

\begin{equation}
\left\{ 
\begin{array}{l}
^{CG}\nabla _{i}^{C}X^{h}=\nabla _{i}X^{h}-\frac{1}{2\alpha }%
R_{.ikm}^{h.}y^{k}X^{m} \\ 
^{CG}\nabla _{\overline{i}}^{C}X^{h}=-\frac{1}{2\alpha }%
R_{.mki}^{h.}y^{k}X^{m} \\ 
^{CG}\nabla _{i}^{C}X^{\overline{h}}=\nabla _{i}\nabla _{k}X^{h}y^{k}-\frac{1%
}{2}R_{imk}^{h}y^{k}X^{m} \\ 
^{CG}\nabla _{\overline{i}}^{C}X^{\overline{h}}=\nabla _{i}X^{h}+\left[ -%
\frac{1}{\alpha }\left( y_{i}\delta _{m}^{h}+y_{m}\delta _{i}^{h}\right) +%
\frac{1+\alpha }{\alpha }g_{im}y^{h}-\frac{1}{\alpha }y_{i}y_{m}y^{h}\right]
X^{h}%
\end{array}%
\right.
\end{equation}

\begin{equation}
^{CG}\nabla _{\beta }^{H}X^{\alpha }=\left( 
\begin{array}{cc}
\nabla _{i}X^{h} & -\frac{1}{2\alpha }R_{.mki}^{h.}y^{k}X^{m} \\ 
-\frac{1}{2}R_{imk}^{h}y^{k}X^{m} & 0%
\end{array}%
\right) .
\end{equation}

Thus, taking account of the fact that $\nabla _{i}X^{h}=0$ implies $%
R_{imk}^{h}X^{m}=0,$ we have

\begin{proposition}
The Complete and horizontal lifts of a vector field in $M_{n}$ to $T\left(
M_{n}\right) $ with the $Cheeger-Gromoll$ metric $^{CG}g$ are parallel if
and only if the given vector field in $M_{n}$ is parallel.
\end{proposition}

We now consider the vertical, complete and horizontal lifts of a vector
field in $M_{n}$ with local components $X^{h}$ in $M_{n}$ to $T\left(
M_{n}\right) $ and compute components of the associated covector fields of $%
X $ \ with respect to the metric $^{CG}g.$Then we obtain respectively

\begin{equation}
\left\{ 
\begin{array}{l}
^{V}X_{B}=\left( 0,\frac{1}{\alpha }\left(
X_{i}+g_{is}X_{t}y^{s}y^{t}\right) \right) \\ 
^{C}X_{B}=\left( X_{i},\frac{1}{\alpha }\left( \nabla X_{i}+g_{is}\nabla
X_{t}y^{s}y^{t}\right) \right) \\ 
^{H}X_{B}=\left( X_{i},0\right)%
\end{array}%
\right.
\end{equation}%
with respect to the adapted frame, where $X_{i}=g_{is}X^{s}$ are local
components of the associated covector field $X^{\ast }$ in $M_{n}.$ Thus we
see that the vertical, complete and horizontal lifts of the associated
covector field $X^{\ast }$ have respectively covariant derivatives with
components

\begin{equation*}
^{CG}\nabla _{\beta }^{V}X_{\gamma }^{{}}=\left( 
\begin{array}{cc}
^{CG}\nabla _{i}^{V}X_{j} & ^{CG}\nabla _{i}^{V}X_{\overline{j}} \\ 
^{CG}\nabla _{\overline{i}}^{V}X_{j} & ^{CG}\nabla _{\overline{i}}^{V}X_{%
\overline{j}}%
\end{array}%
\right)
\end{equation*}

\begin{equation}
\left\{ 
\begin{array}{l}
^{CG}\nabla _{i}^{V}X_{j}=-\frac{1}{2\alpha }R_{ijk}^{h}y^{k}\left(
X_{h}+g_{hs}X_{t}y^{s}y^{t}\right) \\ 
^{CG}\nabla _{i}^{V}X_{\overline{j}}=\frac{1}{\alpha }\left( \nabla
_{i}X_{j}+\nabla iX_{t}g_{js}y^{s}y^{t}\right) \\ 
^{CG}\nabla _{\overline{i}}^{V}X_{j}=0 \\ 
^{CG}\nabla _{\overline{i}}^{V}X_{\overline{j}}=\left[ -\left( y_{i}\delta
_{j}^{h}+y_{j}\delta _{i}^{h}\right) +\left( 1+\alpha \right)
g_{ij}y^{h}-y_{i}y_{j}y^{h}\right] \frac{\left(
X_{h}+g_{hs}X_{t}y^{s}y^{t}\right) }{\alpha ^{2}}%
\end{array}%
\right.
\end{equation}

\begin{equation*}
^{CG}\nabla _{\beta }^{C}X_{\gamma }^{{}}=\left( 
\begin{array}{cc}
^{CG}\nabla _{i}^{C}X_{j} & ^{CG}\nabla _{i}^{C}X_{\overline{j}} \\ 
^{CG}\nabla _{\overline{i}}^{C}X_{j} & ^{CG}\nabla _{\overline{i}}^{C}X_{%
\overline{j}}%
\end{array}%
\right)
\end{equation*}

\bigskip 
\begin{equation}
\left\{ 
\begin{array}{l}
^{CG}\nabla _{i}^{C}X_{j}=\nabla _{i}X_{j}-\frac{1}{2\alpha }%
R_{ijk}^{h}y^{k}\left( \nabla X_{h}+g_{hs}\nabla X_{t}y^{s}y^{t}\right) \\ 
^{CG}\nabla _{i}^{C}X_{\overline{j}}=\frac{1}{\alpha }\left( \nabla
_{i}\nabla _{n}X_{j}y^{n}+\left( g_{js}\nabla _{i}\nabla
_{n}X_{t}y^{s}y^{t}y^{n}\right) \right) -\frac{1}{2\alpha }%
R_{.ikj}^{m.}y^{k}X_{m} \\ 
^{CG}\nabla _{\overline{i}}^{C}X_{j}=-\frac{1}{2\alpha }%
R_{.jki}^{m.}y^{k}X_{m} \\ 
^{CG}\nabla _{\overline{i}}^{C}X_{\overline{j}}=\frac{1}{\alpha }\left(
\nabla _{i}X_{j}+g_{js}\nabla _{i}X_{t}y^{s}y^{t}\right) \\ 
\text{ \ }+\left[ \left( y_{i}\delta _{j}^{h}+y_{j}\delta _{i}^{h}\right)
-\left( 1+\alpha \right) g_{ij}y^{h}+y_{i}y_{j}y^{h}\right] \frac{\left(
\nabla X_{h}+g_{ms}X_{t}y^{s}y^{t}\right) }{\alpha ^{2}}%
\end{array}%
\right.
\end{equation}

\begin{equation}
^{CG}\nabla _{\beta }^{H}X_{\gamma }^{{}}=\left( 
\begin{array}{cc}
\nabla _{i}X_{j} & -\frac{1}{2\alpha }R_{.ikj}^{m.}y^{k}X_{m} \\ 
-\frac{1}{2\alpha }R_{.jki}^{m.}y^{k}X_{m} & 0%
\end{array}%
\right)
\end{equation}

with respect to the adapted frame. Thus the rotations of $^{H}X,^{C}X$ and $%
^{V}X$ have respectively components of the form

\begin{equation}
^{CG}\nabla _{\beta }^{H}X_{\gamma }^{{}}-^{CG}\nabla _{\gamma }^{H}X_{\beta
}^{{}}=\left( 
\begin{array}{cc}
\nabla _{i}X_{j}-\nabla _{j}X_{i} & 0 \\ 
0 & 0%
\end{array}%
\right) ,
\end{equation}

\begin{equation*}
^{CG}\nabla _{\beta }^{C}X_{\gamma }^{{}}-^{CG}\nabla _{\gamma }^{C}X_{\beta
}^{{}}=\left( 
\begin{array}{cc}
\mathit{A} & \mathit{B} \\ 
\mathit{C} & \mathit{D}%
\end{array}%
\right)
\end{equation*}

\begin{equation}
\left\{ 
\begin{array}{l}
\mathit{A=}\left( \nabla _{i}X_{j}-\nabla _{j}X_{i}\right) -\frac{1}{\alpha }%
R_{ijk}^{h}y^{k}\left( \nabla X_{h}+g_{hs}\nabla X_{t}y^{s}y^{t}\right) \\ 
\mathit{B}=\frac{1}{\alpha }\left( \nabla _{i}\nabla _{n}X_{j}y^{n}+\left(
g_{js}\nabla _{i}\nabla _{n}X_{t}-g_{is}\nabla _{j}\nabla _{n}X_{t}\right)
y^{s}y^{t}y^{n}\right) \\ 
\mathit{C}=0 \\ 
\mathit{D}=\frac{1}{\alpha }\left[ \left( \nabla _{i}X_{j}-\nabla
_{j}X_{i}\right) +\left( g_{js}\nabla _{i}X_{t}-g_{is}\nabla
_{j}X_{t}\right) y^{s}y^{t}\right] \\ 
\end{array}%
,\right.
\end{equation}

\begin{eqnarray*}
^{CG}\nabla _{\beta }^{V}X_{\gamma }^{{}}-^{CG}\nabla _{\gamma }^{V}X_{\beta
}^{{}} &=&\left( 
\begin{array}{cc}
\mathit{A}^{\prime } & \mathit{B}^{\prime } \\ 
\mathit{C}^{\prime } & \mathit{D}^{\prime }%
\end{array}%
\right) \\
&&
\end{eqnarray*}

\begin{equation}
\left\{ 
\begin{array}{l}
\mathit{A}^{\prime }\mathit{=}-\frac{1}{\alpha }R_{ijk}^{h}y^{k}\left(
X_{h}+g_{hs}X_{t}y^{s}y^{t}\right) \\ 
\mathit{B}^{\prime }=\frac{1}{\alpha }\left[ \left( \nabla _{i}X_{j}-\nabla
_{j}X_{i}\right) +\left( g_{js}\nabla _{i}X_{t}-g_{is}\nabla
_{j}X_{t}\right) y^{s}y^{t}\right] \\ 
\mathit{C}^{\prime }=0 \\ 
\mathit{D}^{\prime }=0%
\end{array}%
\right.
\end{equation}%
with respect to the adapted frame.

From $\left( 12\right) $, we see, If complete lift of the associated
covector field of $X$ is closed in $T\left( M_{n}\right) ,$then

\begin{equation}
\nabla _{i}X_{j}-\nabla _{j}X_{i}=0,\text{ \ \ }\nabla _{i}\nabla _{j}X_{t}=0%
\text{\ . \ }
\end{equation}%
Further, if the conditions$\left( 14\right) $ are satisfied, we deduce that

\begin{equation*}
R_{ijk}^{h}\left( \nabla X_{h}+g_{hs}\nabla X_{t}\right) =0.
\end{equation*}

From there we have

\begin{proposition}
The complete lift of the associated covector field of $X$ is closed in $%
T\left( M_{n}\right) $ if the associated covector field of $X$ is closed and
the second covariant derivative of $X$ vanishes in $M_{n}.$
\end{proposition}

From $\left( 11\right) $, we see that the horizontal lift of the associated
covector field of $X$ is closed in $T\left( M_{n}\right) $ if and only if
the horizontal lift of the associated covector field of $X$ \ is closed.

We now compute the Lie derivatives of the metric $^{CG}g$ with respect to $%
^{C}X$ and $^{H}X$, by making use of $\left( 7\right) $. The Lie derivatives
of the metric $^{CG}g$ with respect to $^{V}X,^{C}X$ and $^{H}X$ have
respectively components

\begin{equation*}
\tciLaplace _{^{C}X}^{{}}{}^{CG}g=^{CG}\nabla _{\beta }^{C}X_{\gamma
}^{{}}+^{CG}\nabla _{\gamma }^{C}X_{\beta }^{{}}=\left( 
\begin{array}{cc}
\mathit{A}_{1} & \mathit{B}_{1} \\ 
\mathit{C}_{1} & \mathit{D}_{1}%
\end{array}%
\right)
\end{equation*}

\begin{equation}
\left\{ 
\begin{array}{l}
\mathit{A}_{1}\mathit{=}\left( \nabla _{i}X_{j}+\nabla _{j}X_{i}\right) \\ 
\mathit{B}_{1}=\frac{1}{\alpha }\left[ \nabla _{i}\nabla
_{n}X_{j}y^{n}+\left( g_{js}\nabla _{i}\nabla _{n}X_{t}+g_{is}\nabla
_{j}\nabla _{n}X_{t}\right) y^{s}y^{t}y^{n}-R_{.ikj}^{m.}y^{k}X_{m}\right]
\\ 
\mathit{C}_{1}=-\frac{1}{\alpha }R_{.jki}^{m.}y^{k}X_{m} \\ 
\mathit{D}_{1}=\frac{1}{\alpha }\left[ \left( \nabla _{i}X_{j}-\nabla
_{j}X_{i}\right) +\left( g_{js}\nabla _{i}X_{t}-g_{is}\nabla
_{j}X_{t}\right) y^{s}y^{t}\right] \\ 
+\left[ \left( y_{i}\delta _{j}^{m}+y_{j}\delta _{i}^{m}\right) -\left(
1+\alpha \right) g_{ij}y^{m}+y_{i}y_{j}y^{m}\right] \frac{2\left( \nabla
X_{m}+g_{ms}X_{t}y^{s}y^{t}\right) }{\alpha ^{2}} \\ 
\end{array}%
\right.
\end{equation}

\begin{equation}
\tciLaplace _{^{H}X}^{{}}{}^{CG}g=^{CG}\nabla _{\beta }^{H}X_{\gamma
}^{{}}+^{CG}\nabla _{\gamma }^{H}X_{\beta }^{{}}=\left( 
\begin{array}{cc}
\left( \nabla _{i}X_{j}+\nabla _{j}X_{i}\right) & -\frac{1}{\alpha }%
R_{.ikj}^{m.}y^{k}X_{m} \\ 
-\frac{1}{\alpha }R_{.jki}^{m.}y^{k}X_{m} & \mathit{0}%
\end{array}%
\right)
\end{equation}%
with respect to the adapted frame in $T\left( M_{n}\right) .$

Since we have

\begin{equation*}
\nabla _{i}\nabla _{n}X_{j}+\left( g_{js}\nabla _{i}\nabla
_{n}X_{t}+g_{is}\nabla _{j}\nabla _{n}X_{t}\right) -R_{.ikj}^{m.}X_{m}=0
\end{equation*}%
as a consequence of $\tciLaplace _{X}$ $g_{ji}=\nabla _{i}X_{j}+\nabla
_{j}X_{i}=0\cite{7}$, we conclude by means of $\left( 15\right) $ that the
complete lift $^{C}X$ is a Killing vector field in $T\left( M_{n}\right) $
if and only if is a Killing vector field in $M_{n}.$

We next have

\begin{equation*}
R_{.ikj}^{m.}X_{m}=0
\end{equation*}%
as a consequence of the vanishing of the second covariant derivative of $X$.
Conversely, the conditions

\begin{equation*}
\nabla _{i}X_{j}+\nabla _{j}X_{i}=0\text{ \ \ \ and \ \ \ }%
R_{.ikj}^{m.}X_{m}=0\text{\ \ \ }
\end{equation*}%
imply that the second covariant derivative of $X$ \ vanishes. From these
results, we have

\begin{proposition}
Necessary and sufficient conditions in order that $\left( a\right) $ the
complete lift to $T\left( M_{n}\right) $, with metric ${}^{CG}g$, of a
vector field $X$ in $M_{n}$ be a Killing vector field in $T\left(
M_{n}\right) $ is that $X$ is a Killing vector field with vanishing
covariant derivative in $M_{n}$, $\left( b\right) $ the horizontal lift to $%
T\left( M_{n}\right) $, with metric ${}^{CG}g$, of a vector field $X$ in $%
M_{n}$ be a Killing vector field in $T\left( M_{n}\right) $ is that $X$ is a
Killing vector field with vanishing second covariant derivative in $M_{n}.$
\end{proposition}

We now calculate the components of the curvature tensor $^{CG}R$ of $T\left(
M_{n}\right) $ with the metric $Cheeger-Gromoll$ metric $^{CG}g.$ Components
of the curvature tensor $^{CG}R$ with respect to the adapted frame are given
by

\begin{equation}
^{CG}R_{\alpha \beta \gamma }^{\delta }=D_{\alpha }{}^{CG}\Gamma _{\beta
\gamma }^{\delta }-D_{\beta }{}^{CG}\Gamma _{\alpha \gamma }^{\delta
}+^{CG}\Gamma _{\alpha \epsilon }^{\delta }{}^{CG}\Gamma _{\beta \gamma
}^{\epsilon }-^{CG}\Gamma _{\beta \epsilon }^{\delta }{}^{CG}\Gamma _{\alpha
\gamma }^{\epsilon }-\Omega _{\alpha \beta }^{\epsilon }{}^{CG}\Gamma
_{\epsilon \gamma }^{\delta }
\end{equation}%
where $D_{\alpha }$ are local vector fields and $\Omega _{\alpha \beta
}^{\epsilon }$ are components of the non-holonomic object. From $\left(
17\right) $ we have

\begin{equation}
\left\{ 
\begin{array}{l}
^{CG}R_{jik}^{h}=R_{jik}^{h}+\frac{1}{4\alpha }\left(
R_{.imn}^{h.}R_{jkl}^{n}-R_{.jmn}^{h.}R_{ikl}^{n}\right) y^{m}y^{l}-\frac{1}{%
2\alpha }R_{jim}^{n}R_{.ksn}^{h.}y^{m}y^{s} \\ 
^{CG}R_{jik}^{\overline{h}}=\frac{1}{2}\left( \nabla _{j}R_{ikm}^{h}-\nabla
_{i}R_{jkm}^{h}\right) y^{m} \\ 
^{CG}R_{ji\overline{k}}^{\overline{h}}=R_{jik}^{h}+\frac{1}{4\alpha }\left(
R_{inm}^{h}R_{.jlk}^{n.}-R_{.jnm}^{h}R_{.ilk}^{n.}\right) y^{m}y^{l} \\ 
\text{ \ \ \ \ \ \ \ \ \ }-R_{jim}^{h}y^{m}\left[ -\frac{1}{\alpha }\left(
y_{n}\delta _{k}^{h}+y_{k}\delta _{n}^{h}\right) +\frac{1+\alpha }{\alpha }%
g_{nk}y^{y}-\frac{1}{\alpha }y_{n}y_{k}y^{h}\right] \\ 
^{CG}R_{\overline{j}\overline{i}k}^{h}=\frac{1}{4\alpha ^{2}}\left(
R_{.nmj}^{h.}R_{.kli}^{n.}-R_{.nmi}^{h.}R_{.klj}^{n.}\right) y^{m}y^{l} \\ 
^{CG}R_{\overline{j}\overline{i}k}^{\overline{h}}=0,^{CG}R_{\overline{j}%
\overline{i}\overline{k}}^{\overline{h}}=0,^{CG}R_{j\overline{i}\overline{k}%
}^{\overline{h}}=0 \\ 
^{CG}R_{j\overline{i}k}^{h}=-\frac{1}{2\alpha }\left( \nabla
_{j}R_{.kmi}^{h.}\right) y^{m} \\ 
^{CG}R_{j\overline{i}k}^{\overline{h}}=\frac{1}{2}R_{jik}^{h}+\frac{1}{%
4\alpha }R_{jnm}^{h}R_{.kli}^{n.}y^{m}y^{l} \\ 
\text{ \ \ \ \ \ \ \ \ \ }+\frac{1}{2\alpha }\left[ \left( y_{i}\delta
_{n}^{h}+y_{n}\delta _{i}^{h}\right) +\frac{1+\alpha }{\alpha }g_{in}y^{n}-%
\frac{1}{\alpha }y_{n}y_{k}y^{h}\right] R_{jkm}^{n}y^{m} \\ 
\multicolumn{1}{c}{}%
\end{array}%
\right.
\end{equation}

with respect to the adapted frame. Thus we have

\begin{proposition}
The tangent bundle $T\left( M_{n}\right) $ over a Riemannian manifold $M_{n}$
is locally flat with respect to the metric $^{CG}g$ if and only if $M_{n}$
is locally flat.
\end{proposition}

\end{document}